\documentclass[final,onefignum,onetabnum]{siamonline250211}

\usepackage{enumitem}
\setlist[enumerate]{leftmargin=.5in}
\setlist[itemize]{leftmargin=.5in}

\usepackage{amsfonts,amsmath}
\usepackage{dsfont}
\usepackage{comment}

\allowdisplaybreaks

\newtheorem{remark}[theorem]{Remark}

\DeclareMathAlphabet{\mathpzc}{OT1}{pzc}{m}{it}
\newcommand \pd {\mathpzc{d}}
\newcommand \bd {\bar{d}}

\newcommand \eps {\varepsilon}
\newcommand \gam {\gamma}

\newcommand \mC {\mathcal{C}}
\newcommand \mL {\mathcal{L}}
\newcommand \mF {\mathcal{F}}

\newcommand \drm {\mathrm{d}}
\newcommand \e {\mathrm{e}}

\newcommand \E {\mathbb{E}}
\newcommand \F {\mathbb{F}}
\renewcommand \P {\mathbb{P}}
\newcommand \R {\mathbb{R}}
\newcommand \V {\mathbb{V}\mathrm{ar}}

\newcommand \xt {\tilde x}
\newcommand \Vt {\widetilde{V}}

\newcommand{\id}{\mathds{1}}

\allowdisplaybreaks

\begin{document}
	
	\title{Equilibrium Mean-Variance Dividend Rate Strategies\thanks{This version: July 28, 2025. Forthcoming in \emph{SIAM Journal on Financial Mathematics}. \funding{The first and second authors acknowledge the financial support from the Natural Sciences and Engineering Research Council of Canada, grants 05061 and 04958, respectively. The third author thanks the Cecil J. and Ethel M. Nesbitt Professorship for the financial support of her research.}}}	
	
	\author{Jingyi Cao%
		\thanks{Department of Mathematics and Statistics, York University, Canada. Email: jingyic@yorku.ca}
		\and Dongchen Li%
		\thanks{Department of Mathematics and Statistics, York University, Canada. Email: dcli@yorku.ca}
		\and Virginia R. Young%
		\thanks{Department of Mathematics, University of Michigan, USA. Email: vryoung@umich.edu}
		\and Bin Zou%
		\thanks{Department of Mathematics, University of Connecticut, USA. Email: bin.zou@uconn.edu}
	}

	\maketitle
	
	\begin{abstract}
		This paper studies an optimal dividend problem for a company that aims to maximize the mean-variance (MV) objective of the accumulated discounted dividend payments up to its ruin time. The MV objective involves an integral form over a random horizon that depends endogenously on the company's dividend strategy, and these features lead to a novel time-inconsistent control problem. To address the time inconsistency, we seek a time-consistent equilibrium dividend rate strategy. We first develop and prove a new verification lemma that characterizes the value function and equilibrium strategy by an extended Hamilton-Jacobi-Bellman system. Next, we apply the verification lemma to obtain the equilibrium strategy and show that it is a barrier strategy for small levels of risk aversion. 
	\end{abstract}
	
	

	\section{Introduction}
	
	Dividend decisions are critical to the operation of a company because they affect the earnings that a company distributes to its shareholders versus the amount it retains and reinvests (Baker \cite{B2009}, p.3). The study of optimal dividend policies dates back to the seminal work of De Finetti \cite{D1957} and has been an active research topic in actuarial science and financial mathematics for decades (see Albrecher and Thonhauser \cite{AT2009} and Avanzi \cite{A2009} for survey articles). This paper formulates a dynamic dividend problem in continuous time under the MV criterion and seeks a time-consistent equilibrium solution.
	
	We consider a company (for instance, an insurer) that pays dividends to its shareholders.  Assume that the company's uncontrolled surplus (that is, excluding dividend payments), $X^0 = (X^0_t)_{t \ge 0}$, is modeled by a Brownian motion with drift, $\drm X_t^0 = a \, \drm t + b \, \drm B_t$.\footnote{Such a model is called the \emph{diffusion approximation model} in risk theory (see Grandell \cite{G1991}), and it is frequently used in the study of optimal dividend problems (see, for instance,  Asmussen and Taksar \cite{AT1997}).}  
	Here, $a \in \R$ and $b > 0$ are constants, and $B = (B_t)_{t \ge 0}$ is a standard, one-dimensional Brownian motion defined on its natural filtration space  $(\Omega, \mF, \F = (\mF_t)_{t \ge 0}, \P)$. 
	We adopt the classical control framework and assume that the company pays dividends at a bounded rate (see Case A in Jeanblanc-Picqu\'e and Shiryaev \cite{JS1995} or Section 2 in Asmussen and Taksar \cite{AT1997}).   
	We restrict our attention to Markov feedback controls induced by \emph{deterministic} functions $\pd: (x, t) \in \R_+^2 \to \R_+$ (see Bj\"ork and Murgoci \cite{BM2010}); that is, the company pays dividend at rate $\pd(x, t)$ at time $t$ when its surplus is $x$. 
	As such, for a given dividend strategy $\pd$, the company's controlled surplus $X = (X_t)_{t \ge 0}$ follows the dynamics
	\begin{align}
		\label{eq:X}
		\drm X_t = \left( a - \pd(X_t, t) \right) \, \drm t + b \, \drm B_t, 
	\end{align}
	with $X_0 = x_0 > 0$.  Define the ruin time of the company by $\tau := \inf \{ t \ge 0: X_t < 0 \}$.  Note that $X$ and $\tau$ depend on the dividend strategy $\pd$ (that is, $X:=X^\pd$ and $\tau := \tau^\pd$), but we omit that dependence for notational simplicity and will follow the same rule hereafter.

	Let $Y_t$ denote the discounted dividends paid between time $t$ and ruin time $\tau$ under strategy $\pd$, using a discount rate $\rho > 0$, that is,
	\begin{equation}\label{eq:Y}
		Y_t = \int_t^\tau \e^{-\rho(s - t)} \, \pd(X_s, s) \, \drm s.
	\end{equation}
	We set $Y_t = 0$ if $t \ge \tau$ (ruin has occurred before or at time $t$). 
	Arguably speaking, the most common objective in the study of optimal dividend is to maximize the (conditional) expectation of $Y_t$; see, for instance, Equation (2.4) in Taksar \cite{T2000}.
	However, as pointed out in Avanzi \cite{A2009} (p.239), ``In the real world, variability in dividend payments (especially, decreases) is not well received in the markets and may penalize significantly the value of the share.''
	Since variance is a popular measure of variability, we incorporate a variance term of $Y_t$ to penalize variability in dividend payments and propose the following MV objective: 
	\begin{equation}\label{eq:J}
		J(x, t; \pd) = \E_{x, t}(Y_t) - \dfrac{\gam}{2} \, \V_{x, t}(Y_t),
	\end{equation}
	in which $\E_{x, t}$ and $\V_{x, t}$ denote expectation and variance under $\P$, respectively, conditional on $X_t = x \ge 0$ and $\tau > t$, and the parameter $\gam > 0$ measures the company's risk aversion toward variance.\footnote{Avanzi et al.\ \cite{AFS2023} also aim to ensure stability in dividend payments, a motivation shared with this paper, but they take a different approach by restricting to affine dividend strategies and adopting a linear-quadratic criterion (with target dividend rates).
		We remark that there is a strand of literature imposing drawdown or ratcheting constraints on dividend strategies that penalize or prevent the decrease in dividends; see, for instance, Angoshtari et al.\ \cite{ABY2019}, Albrecher et al.\ \cite{AAM2022}, and Guan and Xu \cite{GX2024}.}  
	Note that when $\gam = 0$,  the dynamic problem of $\max \, \E_{x, t}(Y_t)$ is time-consistent,  and can be solved by the standard Hamilton-Jacobi-Bellman (HJB) method; see Theorem 2.1 in Taksar \cite{T2000}.
	However, maximizing the dynamic MV objective in \eqref{eq:J} with $\gam > 0$ for all $(x, t) \in \R_+^2$ leads to a \emph{time-inconsistent} problem because variance does not satisfy iterated expectation, and in turn, the Dynamic Programming Principle fails (see Basak and Chabakauri \cite{BC2010} and Bj\"ork and Murgoci \cite{BM2010} for a detailed discussion). To address the time-inconsistency issue, we follow the game-theoretic approach, as in Bj\"ork and Murgoci \cite{BM2010}, in which the company manager plays against all future versions of themselves and achieves consistent planning via a Nash game. We seek a time-consistent equilibrium strategy $\pd^*$ (see \cref{def:bq_timeconsist}); this problem is new to the optimal-dividend literature.
	
	The objective $J$ in \eqref{eq:J} differs from the standard MV objective in two noticeable ways.  Specifically, (1) $Y_t$ in \eqref{eq:Y} is an integral of random variables, and (2) it is defined up to the ruin time $\tau$, which depends \emph{endogenously} on the company's dividend strategy $\pd$; by comparison, the standard MV objective considers one single random variable at the terminal (constant) time $T > 0$ (see, for instance, Bj\"ork et al.\ \cite{BMZ2014}). Because of these differences, the pursuit of an equilibrium strategy $\pd^*$ in this paper turns out to be challenging; in fact, our problem in  \cref{def:bq_timeconsist} is \emph{not} covered even by the most general case in Bj\"ork and Murgoci \cite{BM2010} (compare $J$ in \eqref{eq:J} with their $J$ in Equation (39)). Therefore, we develop and prove a verification theorem that is tailored for our MV dividend problem; see \cref{thm:verif_bd}. 
	To the best of our knowledge, this verification theorem is new to the time-inconsistent control literature.  
	In particular, the extended HJB system in \eqref{eq:HJB}-\eqref{eq:H} involves three functions, $V(x,t) = J(x,t; \pd^*)$, $G=\E_{x, t}(Y_t^*)$, and $H = \E_{x, t} ((Y_t^*)^2 )$, in which $Y^* := Y^{\pd^*}$; note that the presence of $H$ in the extended HJB system is due to the integral form of $Y_t$ in \eqref{eq:Y} (see Step 4 in the proof of Theorem 2.3).
	On the technical side, several difficulties arise in the proof of \cref{thm:verif_bd}. First, we need to work out a second-order PDE for $H$ in \eqref{eq:H} and an HJB equation for $\Vt$ in \eqref{eq:HJB} and show that their solutions equal $\E_{x, t} ((Y_t^*)^2 )$ and $V$ (the equilibrium value function), respectively. 
	Second, it takes a delicate analysis to obtain the first-order expansion of the objective under perturbed strategies and to show that the optimizer to \eqref{eq:HJB} is an equilibrium strategy (Step 4 in the proof). 
	Although Landriaut et al.\ \cite{LLLY2018} also seek equilibrium strategies to an MV optimization problem under a \emph{random} horizon, they consider the controlled process only at the random time $\tau$ (that is, $X_\tau$), and $\tau$ in their paper is \emph{exogenously} given and \emph{independent} of the driving Brownian motion.  
	Kronborg and Steffensen \cite{KS2015} formulate a class of time-inconsistent MV  investment-consumption problems, in which the MV criterion is applied to the discounted wealth at the terminal time $T$ plus an integral of discounted intertemporal consumption over $[t, T]$. They allow the decision maker to continue investing and consuming after ruin; by contrast, our model ``ends the game'' if ruin occurs, as in the usual dividend problem. Since they solve a finite-horizon problem, the value function and the equilibrium strategies, in general, depend on both $t$ and $T$ (see their Proposition 4.1). But our random-horizon problem is similar to an infinite-horizon problem, in the sense that the value function and the equilibrium (barrier) strategy are time-invariant (see Theorems \ref{thm:small_gam_bd} and \ref{thm:max_d}).
	
	Next, we apply the verification theorem (\cref{thm:verif_bd}) to obtain the company's equilibrium dividend strategy $\pd^*$. We show that for small risk aversion $\gam$, $\pd^*$ is a barrier strategy; namely, the company pays at the maximum rate when its surplus is above the barrier $\xt$ and pays nothing when it is below $\xt$ (see \cref{thm:small_gam_bd}). Given $\gam$, $\xt$ is characterized as the unique solution to a non-linear equation $f(x, \gam) - 1 = 0$ (see \cref{lem:small}) and can be found efficiently by a standard non-linear solver (see \cref{fig:gamma}). Furthermore, we identify a sufficient condition under which the company always pays dividends at the maximum rate (see \cref{thm:max_d}). 
	
	As hinted earlier, we are not aware of any study in optimal dividend problems that seeks equilibrium strategies under time-inconsistent MV preferences as ours in \eqref{eq:J}. But time inconsistency might also arise from sources that are different from the variance term in MV preferences, and a prime example is non-exponential discounting. Several papers solve for equilibrium dividend strategies under quasi-hyperbolic (or equivalently piecewise exponential) discount functions; see, for instance, Chen et al.\ \cite{CLZ2014, CLZ2018}, Zhu et al.\ \cite{ZSY2020}, and Zhou and Jin \cite{ZJ2020}.

	The rest of the paper is organized as follows. We prove a verification theorem for the company's MV dividend problem in \cref{sec:veri} and apply it to obtain the company's equilibrium dividend strategy in \cref{sec:eq}. Finally, \cref{sec:con} concludes the study.

	\section{Verification theorem}
	\label{sec:veri}
	
	We first introduce the set of admissible strategies and formally define time-consistent equilibrium strategies below.
	
	\begin{definition}\label{def:pd}
		A dividend strategy $\pd = (\pd(X_t, t))_{t \ge 0}$ is called {\rm admissible} if $\pd$ is a deterministic, Borel-measurable function taking values in $[0, \bd]$, for some maximum rate $\bd > 0$.
	\end{definition}
	
	\begin{definition}\label{def:bq_timeconsist} 
		Let $\pd^* = \left(\pd^*(X_s^*, s)\right)_{s \ge 0}$ be an admissible strategy, in which $X^* := X^{\pd^*}$ is the company's surplus under strategy $\pd^*$.  
		Fix an arbitrary initial time $t \ge 0$, $\eps>0$, and $d \in [0, \bd]$, and define the perturbed strategy $\pd^\eps = \left( \pd^\eps(X_s, s) \right)_{s \ge t}$ by
		\begin{equation}\label{eq:d_eps}
			\pd^\eps(X_s, s) =
			\begin{cases}
				d, &\quad t \le s < (t + \eps) \wedge \tau, \\
				\pd^*(X_s, s), &\quad s \ge (t + \eps) \wedge \tau,
			\end{cases}
		\end{equation}
		in which $\tau := \tau^{\pd^\eps}$ and $X := X^{\pd^\eps}$.
		The strategy $\pd^*$ is said to be a {\rm time-consistent equilibrium dividend strategy} if, for all $(x, t) \in \R_+^2$,
		\begin{equation}\label{eq:J_eqm}
			\liminf_{\eps \to 0^+} \, \dfrac{J(x, t; \pd^*) - J(x, t; \pd^\eps)}{\eps} \ge 0.
		\end{equation}
		If an equilibrium strategy $\pd^*$ exists, we call $V(x, t) = J(x, t; \pd^*)$ the {\rm equilibrium value function}.
	\end{definition}
	
	Assume for a moment that starting from time $t+\eps$ onward, the company manager will follow $\pd^*$ to pay dividends. If such a strategy $\pd^*$ is an equilibrium strategy, then the manager should have \emph{no} incentive to deviate from it, which is captured by the first-order condition in (2.2).  Indeed, (2.2) implies that  $J(x, t; \pd^*) \ge J(x, t; \pd^\eps) + o(\eps)$ for all ``deviation'' $\pd^\eps$, and, thus, the manager should ``stick to'' $\pd^*$ over $[t, t+\eps)$ as well, making $\pd^*$ a consistent planning over the entire time horizon.

	Before we state the verification theorem for the MV dividend problem in \cref{def:bq_timeconsist}, we define a differential operator $\mL^{d}$, for all $\phi \in \mC^{2, 1}(\R_+^2)$ and $d \in [0, \bd]$, by
	\begin{align*}
		\mL^{d} \, \phi(x, t) &= \partial_t \phi(x, t)  +  (a - d) \partial_x \phi(x, t) + \dfrac{1}{2} \, b^2  \partial_{xx} \phi(x, t),
	\end{align*}
	in which $\partial_t$, $\partial_x$, and $\partial_{xx}$  denote the corresponding partial derivatives of $\phi$.

	\begin{theorem}\label{thm:verif_bd}	
		Suppose there exist three functions $\Vt$, $G$, and $H$, all mapping from $(x, t) \in \R_+^2 \to \R$, that satisfy the following conditions: 
		\begin{enumerate}
			\item $\Vt, G , H \in \mC^{2, 1}(\R_+^2)$, except that $G(\cdot, t)$ and $H(\cdot, t)$ might only be $\mC^1$ along a specific curve $x = \xt(t)$ for all $t \ge 0$ with left and right second derivatives. 
			$\Vt$, $G$, and $H$ satisfy regularity conditions such that the stochastic integrals in \eqref{eq:G_rando} and \eqref{eq:hatH_rando} are martingales and $\lim_{s \to \infty} \E_{x, t} \left( \e^{-\rho s } \, \phi(X_s, s) \right) = 0$ for $\phi = G, H$.

			\item For all $(x, t) \in \R_+^2$, $\Vt$, $G$, and $H$ jointly solve the following extended HJB system: 
			\begin{align}
				\sup \limits_{d \in [0, \bd]} \big\{ \mL^{d} \, \Vt(x, t) - \frac{\gam}{2} \, \mL^d \, G^2(x, t) + \gam G(x, t) \mL^d \, G(x, t) &  \notag \\
				\; + \, d - \rho G(x, t) + \gam \rho \big(H(x, t) - G^2(x, t) \big) \big\}  &= 0,   \label{eq:HJB}  \\
				\mL^{\pd^*(x, t)} G(x, t) - \rho G(x, t) + \pd^*(x, t) &= 0,  \label{eq:G} \\
				\mL^{\pd^*(x, t)} H(x, t) - 2 \rho H(x, t) + 2 \pd^*(x, t) G(x, t) &= 0, \label{eq:H}
			\end{align}
			with boundary conditions $\Vt(0, t) = G(0, t) = H(0, t) = 0$ for all $t \in \R_+$. In \eqref{eq:G} and \eqref{eq:H}, $\pd^*(x, t)$ denotes the maximizer of \eqref{eq:HJB} for every $(x, t) \in \R_+^2$.
		\end{enumerate}
		
		\noindent
		Define a dividend strategy $\pd^*$ by $ \left( \pd^*(X_s^*, s) \right)_{s \ge 0}$, in which $X^* = (X_s^*)_{s \ge 0}$ is the company's surplus in \eqref{eq:X} under $\pd^*$, and assume that $\pd^*$ is admissible.  Then, $\pd^*$ is a time-consistent equilibrium dividend strategy, and $\Vt$ equals the company's equilibrium value function $V$  (\cref{def:bq_timeconsist}).
		Moreover, $G$ and $H$ have the following probabilistic representations:
		\begin{align}
			\label{eq:GH}
			G(x, t) = \E_{x, t}(Y_t^*) \quad \text{and} \quad  H(x, t) = \E_{x, t} \left( (Y_t^*)^2 \right),  
		\end{align}
		in which $Y_t^*$ is defined by \eqref{eq:Y} under $\pd^*$; thus, $V(x, t) = G(x, t) - \dfrac{\gam}{2} \left( H(x, t) - G^2(x, t) \right)$. 
	\end{theorem}
	
	\begin{proof}
		Suppose that $\Vt$, $G$, and $H$ satisfy the stated conditions in the theorem, and an equilibrium strategy $\pd^*$ exists. We prove the results in four steps.
		
		\medskip 
		
		\textbf{Step 1.} In this step, we show that if $G$ solves \eqref{eq:G} with $G(0, t) = 0$, then $G(x, t) = \E_{x, t}(Y_t^*)$ in \eqref{eq:GH}.
		The proof is standard in control theory (see, for example, Fleming and Soner \cite{FS2006}), but we include it because we want to refer to parts of it in Step 2 below.  
		
		Let $k > t$ be a fixed number; then, by applying It\^o's formula to $\e^{-\rho(\cdot - t)} G(X^*_{\cdot}, \cdot)$ and using \eqref{eq:G}, we obtain
		\small
		\begin{align}
			\e^{-\rho((\tau \wedge k) - t)} G(X^*_{\tau \wedge k}, \tau \wedge k)
			&= G(X^*_{\tau \wedge t}, {\tau \wedge t}) - \int_{\tau \wedge t}^{\tau \wedge k} \e^{-\rho (s - t)} \pd^* \, \drm s  +  \int_{\tau \wedge t}^{\tau \wedge k} \e^{-\rho (s - t)} b \, \partial_x G \, \drm B_s. \label{eq:G_rando}
		\end{align}
		\normalsize
		Taking conditional expectation, given $X_t^* = x$ and $\tau > t$, on both sides of \eqref{eq:G_rando} yields
		$G(x, t) = \E_{x, t} \int_t^{\tau \wedge k} \e^{- \rho(s - t)} \, \pd^*(X_s^*, s) \, \drm s  + \E_{x, t} \big( \e^{-\rho(k - t)} G(X^*_{k}, k) \id_{\{ \tau > k \}} \big)
		= \E_{x, t} \big( \int_t^{\tau} \e^{- \rho(s - t)} \, \pd^*(X_s^*, s) \, \drm s \big) \allowbreak = \E_{x, t}(Y_t^*),$
		in which the first equality uses $G(X^*_\tau, \tau) = G(0, \tau) = 0$ on $\{\tau \le  k\}$, and the second equality follows by applying the monotone convergence theorem when letting $k \to \infty$ and by using the growth condition of $G$ in Condition 1 to claim that the second limit converges to 0.
		
		\medskip
		
		\textbf{Step 2.} In this step, we show that if $H$ solves \eqref{eq:H} with $H(0,t)=0$, then $	H(x, t) = \E_{x, t} (Y_t^*)^2 $ as claimed in \eqref{eq:GH}.
		Defining $\hat H(x, t) = \e^{-2\rho t} H(x, t)$ and $\hat G(x, t) = \e^{-\rho t} G(x, t)$, 
		we follow an argument similar to the one in Step 1 to obtain 
		\begin{equation}\label{eq:hatG_rando}
			\hat G(X_s^*, s) =  \int_s^{\tau \wedge k} \e^{-\rho u} \pd^*(X_u^*, u) \, \drm u -  \int_s^{\tau \wedge k}  b \, \partial_x \hat{G} (X_u^*, u) \, \drm B_u + \hat G(X^*_{\tau \wedge k}, \tau \wedge k).
		\end{equation}
		Also, \eqref{eq:H} implies that $\hat H$ 
		satisfies $\hat H(0, t) = 0$ and 
		$\mL^{\pd^*(x, t)} \hat H(x, t) + 2\pd^*(x, t) \e^{- \rho t} \hat G(x, t) = 0$.

		Let $k > t$ be a fixed number; then, by It\^o's formula and using the above equality, we obtain 
		\small 
		\begin{align}
			\hat H(X^*_{\tau \wedge k}, \tau \wedge k) 
			= \hat H(X_{t}^*, {t}) - 2 \int_{t}^{\tau \wedge k} \e^{-\rho s} \pd^*(X_s^*, s) \, \hat G(X_s^*, s) \, \drm s   +  \int_{t}^{\tau \wedge k} b \, \partial_x \hat H(X_s^*, s) \, \drm B_s.
			\label{eq:hatH_rando}
		\end{align}
		\normalsize
		Taking conditional expectation on both sides of \eqref{eq:hatH_rando} yields 
		$\hat H(x, t) = \E_{x, t} \big( \hat H(X^*_{\tau \wedge k}, \tau \wedge k) \big)  + 2 \, \E_{x, t} \big(\int_t^{\tau \wedge k} \e^{-\rho s} \pd^*(X_s^*, s) \, \hat G(X_s^*, s) \, \drm s  \big) = \E_{x, t} \big( \hat H(X^*_{\tau \wedge k}, \tau \wedge k) \big) 
		+
		2 \, \E_{x, t} \big( \int_t^{\tau \wedge k} \e^{-\rho s} \pd^*(X_s^*, s) \cdot \allowbreak \big(  \int_s^{\tau \wedge k} \e^{-\rho u} \pd^*(X_u^*, u) \, \drm u \big) \drm s \big) 
		+ 2 \, \E_{x, t} \big( \hat G(X^*_{\tau \wedge k}, \tau \wedge k) \int_t^{\tau \wedge k} \e^{-\rho s} \pd^*(X_s^*, s) \, \drm s \big)$,
		%
		in which the second equality follows by replacing $\hat{G}$ with the right-hand side of  \eqref{eq:hatG_rando} and by using Condition 1 to claim that the related stochastic integral has zero expectation. Then, by letting $k \to \infty$ and by using the growth condition of $G$ and $H$, we obtain 
		\\
		$\hat H(x, t) = 2 \, \E_{x, t} \big( \int_t^{\tau} \e^{-\rho s} \pd^*(X^*_s, s) \left(  \int_s^{\tau} \e^{-\rho u} \pd^*(X_u^*, u) \, \drm u \right) \drm s \big) = \E_{x, t} \left( (Y_t^*)^2 \right) ,$
		in which the second equality follows from the elementary result that 
		\begin{equation}\label{eq:Y_sq}
			2 \int_t^\tau \e^{-\rho s} \pd(X_s, s) \Big( \int_s^\tau \e^{-\rho u} \pd(X_u, u) \, \drm u \Big) \drm s = \Big(\int_t^\tau \e^{- \rho s} \, \pd(X_s, s) \, \drm s \Big)^2
		\end{equation}
		holds for any admissible dividend strategy $\pd$.  
		
		\medskip 
		\textbf{Step 3.}  In this step, we show that if $\Vt$ satisfies \eqref{eq:HJB} with $\Vt(0,t) = 0$, then $\Vt(x, t) = J(x, t; \pd^*)$ for all $(x, t) \in \R_+^2$. By using \eqref{eq:HJB}, \eqref{eq:G}, and \eqref{eq:H}, we obtain 
		\begin{align*}
			\mL^{\pd^*(x, t)}  \Big( \Vt(x, t) - \dfrac{\gam}{2} \, G^2(x, t) \Big) = \mL^{\pd^*(x, t)}  \Big( G(x, t) - \dfrac{\gam}{2} \, H(x, t) \Big).
		\end{align*}
		Thus, by noting the boundary condition $\Vt(X^*_\tau, \tau) - \frac{\gam}{2} \, G^2(X^*_\tau, \tau) = 0 = G(X^*_\tau, \tau) - \frac{\gam}{2} H(X^*_\tau, \tau)$, 
		we deduce $\Vt(x, t)  = G(x, t) - \frac{\gam}{2} \, (H(x, t) - G^2(x, t)) = J(x, t; \pd^*)$.

		\medskip
		
		\textbf{Step 4.} It remains to show that $\pd^*$ is an equilibrium strategy.  To that end, define the strategy $\pd^\eps$ as in \eqref{eq:d_eps}, and we want to prove that the limit in \eqref{eq:J_eqm} holds. For notational easiness, we write $X := X^{\pd^\eps}$ and $Y := Y^{\pd^\eps}$ in this step.
		
		First, by using the definition of $J$ in \eqref{eq:J} and \eqref{eq:Y_sq}, we calculate the objective value under strategy $\pd^\eps$ by (denoting $\pd^\eps_s := \pd^\eps(X_s, s)$)
		\begin{align*}
			J(x, t; \pd^\eps) &= \E_{x, t} \left( \int_t^{\tau} \e^{-\rho (s - t)} \pd^\eps_s \,\drm s \right) - \gam \E_{x, t} \left( \int_t^\tau \e^{-\rho (s - t)} \pd^\eps_s \bigg( \int_s^\tau \e^{-\rho (u - t)} \pd^\eps_s \, \drm u \bigg) \drm s \right) \\
			&\quad + \frac{\gam}{2} \left(\E_{x, t} \left(\int_t^{\tau} \e^{-\rho (s - t)} \pd^\eps_s \, \drm s  \right) \right)^2.
		\end{align*}
		In what follows, we consider each of the three terms above separately and expand them to order $o(\eps)$.  The first term becomes
		\small 
		\begin{align*}
			&\E_{x, t}\Big(\int_t^{\tau} \e^{-\rho (s - t)} \pd^\eps(X_s, s) \drm s \Big) \notag \\
			&= \E_{x, t}\bigg(\int_t^{t + \eps} \e^{-\rho (s - t)} d \, \id_{\{ \tau > s\}} \drm s + \e^{-\rho \eps} \id_{\{ \tau > t + \eps\}}  \, \E_{X_{t+\eps}, t+\eps}\bigg(\int_{t+\eps}^{\tau} \e^{-\rho (s - (t + \eps))}\pd^*(X_s, s) \drm s  \bigg) \bigg) \notag \\
			&= \E_{x, t}\big( \eps d + (1 - \rho \eps) \id_{\{ \tau > t + \eps\}} G(X_{t+\eps}, t+\eps) \big) + o(\eps) \notag \\
			&= \eps d + (1 - \rho \eps) \bigg( \P_{x, t}(\tau > t + \eps) G(x, t) + \E_{x, t} \bigg(  \id_{\{ \tau > t + \eps\}} \int_t^{t+\eps} \mL^d G(X_s, s) \drm s \bigg) \bigg) + o(\eps) \notag \\
			&= \eps d + (1 - \rho \eps) \P_{x, t}(\tau > t + \eps) \big( G(x, t) + \eps \mL^d G(x, t) \big) + o(\eps) \notag \\
			&= \eps d + \P_{x, t}(\tau > t + \eps) \big( G(x, t) + \eps\big(\mL^d G(x, t) - \rho G(x, t) \big) \big) + o(\eps) \notag \\
			&= G(x, t) + \eps \big(\mL^d G(x, t) - \rho G(x, t) + d \big) + o(\eps),
		\end{align*}
		\normalsize 
		in which the last line follows from the Appendix in Grandell \cite{G1991} on finite-time ruin probabilities, specifically,
		\begin{equation}\label{eq:ruin_prob}
			\P_{x, t}(\tau > t + \eps) \sim 1 - \dfrac{b\sqrt{\eps}}{x} \, \exp\bigg( - \dfrac{1}{2\eps} \left( \dfrac{x}{b} \right)^2 \bigg) = 1 + o(\eps).
		\end{equation}
		By using \eqref{eq:ruin_prob} to justify omitting $\id_{\{ \tau > t + \eps\}}$ in the following derivation for simplicity, the second term becomes (ignoring the factor of $-\gam$ for now)
		\small
		\begin{align*}
			& \E_{x, t} \int_t^\tau \e^{-\rho (s - t)} \pd^\eps(X_s, s) \bigg(\int_s^\tau \e^{-\rho (u - t)} \pd^\eps(X_u, u) \, \drm u \bigg) \drm s  \notag \\
			&= \E_{x, t} \int_t^{t + \eps} \e^{-\rho (s - t)} d \bigg( \int_s^{t+\eps} \e^{-\rho (u - t)} d \, \drm u \bigg) \drm s  + \E_{x, t} \int_t^{t + \eps} \e^{-\rho (s - t)} d \bigg( \int_{t+\eps}^\tau \e^{-\rho (u - t)} \pd^*(X_u, u) \, \drm u \bigg) \drm s  \notag \\
			&\quad + \E_{x, t} \int_{t + \eps}^\tau \e^{-\rho (s - t)} \pd^*(X_s, s) \int_s^\tau \e^{-\rho (u - t)} \pd^*(X_u, u) \, \drm u \, \drm s  \notag \\
			&= d^2 \e^{-2\rho \eps} \, \E_{x, t} \int_t^{t + \eps} \e^{-\rho (s - (t+\eps))} \bigg( \int_s^{t+\eps} \e^{-\rho (u - (t+\eps))}  \drm u \bigg) \drm s  \notag \\
			&\quad + d \, \e^{-\rho \eps} \E_{x, t} \int_t^{t + \eps} \e^{-\rho (s - t)} \bigg( \int_{t+\eps}^\tau \e^{-\rho (u - (t+\eps))} \pd^*(X_u, u) \, \drm u \bigg) \drm s \notag \\
			&\quad + \e^{-2\rho \eps} \E_{x, t} \int_{t + \eps}^\tau \e^{-\rho (s - (t + \eps))} \pd^*(X_s, s) \bigg( \int_s^\tau \e^{-\rho (u - (t + \eps))} \pd^*(X_u, u) \, \drm u \bigg) \drm s  \notag \\
			&=  \eps^2 d^2 (1 - 2 \rho \eps)/2 + \eps d (1 - \rho \eps) \E_{x, t}( G(X_{t+\eps}, t+\eps) ) + (1 - 2\rho \eps)/2 \, \E_{x, t}( H(X_{t+\eps}, t+\eps) ) + o(\eps) \notag \\
			&= \eps d \big( G(x, t) + \eps \mL^d G(x, t) \big) + \left(1/2 - \rho \eps \right) \, \big( H(x, t) + \eps \mL^d H(x, t) \big) + o(\eps) \notag \\
			&= H(x, t)/2 + \eps/2 \big(  \mL^d H(x, t) - 2 \rho H(x, t) + 2dG(x, t) \big) + o(\eps).
		\end{align*}
		\normalsize
		By using the result on the first term, the third term becomes (ignoring $\frac{\gam}{2}$ for now)
		\begin{align*}
			\Big(\E_{x, t} \Big(\int_t^{\tau} \e^{-\rho (s - t)} d^\eps(X_s, s) \drm s  \Big) \Big)^2 
			&= \left( G(x, t) + \eps \big(\mL^d G(x, t) - \rho G(x, t) + d \big) \right)^2 + o(\eps) \notag \\
			&= G^2(x, t) + 2 \eps G(x, t) \big(\mL^d G(x, t) - \rho G(x, t) + d \big) + o(\eps).
		\end{align*}
		With the above expansion results, $J(x, t; \pd^\eps)$ becomes
		\begin{align*}
			J(x, t; \pd^\eps) 
			&= \big( G(x, t) - \gam/2 \left( H(x, t) - G^2(x, t) \right) \big) + \eps \big( \mL^d G(x, t) - \rho G(x, t) + d \big) \notag \\
			&\quad + \eps \gam G(x, t) \big( \mL^d G(x, t) - \rho G(x, t) \big) -  \eps \, \gam/2 \, \big( \mL^d H(x, t) - 2 \rho H(x, t) \big) + o(\eps). \hspace{1.5em}
		\end{align*}
		Next, substituting $\frac{\gam}{2} \, H(x, t) = G(x, t) - \Vt(x, t) + \frac{\gam}{2} \, G^2(x, t)$ from Step 3 into the above $J(x, t; \pd^\eps)$, we obtain 
		\begin{align*}
			J(x, t; \pd^\eps) 
			&= \Vt(x, t) + \eps \big( \mL^{d} \, \Vt(x, t) - \gam/2 \, \mL^d \, G^2(x, t) + \gam G(x, t) \mL^d \, G(x, t) \\ 
			&\qquad\qquad \qquad \;\;  + \, d - \rho G(x, t) + \gam \rho \big(H(x, t) - G^2(x, t) \big)  \big) + o(\eps) \\
			&\le \Vt(x, t) + o(\eps) = J(x, t; \pd^*) + \eps , 
		\end{align*}
		in which the inequality follows from \eqref{eq:HJB}.  We have thereby proved the limit in \eqref{eq:J_eqm}.
	\end{proof}

	\begin{remark}
		\label{rem:c1}
		In Condition 1 of \cref{thm:verif_bd}, $G(\cdot, t)$ and $H(\cdot, t)$ might only be $\mC^1$ along a specific curve $x = \xt(t)$. But this poses no issues because $\mL^{\pd(x, t)} G(x, t)$ and $\mL^{\pd(x, t)} H(x, t)$ are piecewise continuous for every admissible strategy $\pd$. 
	\end{remark}

	\section{Equilibrium dividend strategy} 
	\label{sec:eq}
	
	In this section, we use \cref{thm:verif_bd} to obtain the equilibrium dividend strategy $\pd^*$.
	For later convenience, define the following constants:
	\begin{align}
		r_1, r_2 &= \dfrac{1}{b^2} \left[  - a \pm \sqrt{a^2 + 2 \rho b^2} \right],  & 
		r_5, r_6 &= \dfrac{1}{b^2} \left[  - (a - \bd) \pm \sqrt{(a - \bd)^2 + 2 \rho b^2} \right], 
		\label{eq:r_one}
		\\
		r_3, r_4 &= \dfrac{1}{b^2} \left[  - a \pm  \sqrt{a^2 + 4 \rho b^2} \right],
		&
		r_7, r_8 &= \dfrac{1}{b^2} \left[  - (a - \bd) \pm \sqrt{(a - \bd)^2 + 4 \rho b^2} \right], 
		\label{eq:r_two}
	\end{align}
	with $r_1,r_3,r_5,r_7 > 0$ being the positive roots and $r_2,r_4,r_6,r_8 < 0$ the negative roots. 
	
	To build intuition, we first consider the case of $\gam = 0$. Note that the corresponding problem is to maximize $\E_{x, t} (Y_t)$, a time-consistent problem that is already solved in the literature. When $\gam = 0$, Theorem 2.1 in Taksar \cite{T2000} shows that the optimal dividend strategy 
	is to pay at the maximum rate $\bd$ when $x > \xt_0$, and to pay nothing when $x \le \xt_0$, for some barrier $\xt_0 \ge 0$. 
	In this case, the barrier $\xt_0$ is strictly positive if and only if 
	\begin{equation}\label{eq:xt_pos}
		\frac{\bd}{\rho} + \frac{1}{r_6} > 0,
	\end{equation}
	in which $r_6 < 0$ is defined by \eqref{eq:r_one}, and $\rho >0$ is the discount rate in \eqref{eq:Y}.  
	
	When $\gam$ is small, we hypothesize that the equilibrium dividend strategy is also a barrier strategy for some barrier $\xt$. First, we present a technical lemma on ``small'' $\gam$.
	
	\begin{lemma}
		\label{lem:small}
		Define a function $f: (x, \gam) \in \R^2 \to \R$ by 
		\begin{align*}
			f(x, \gam) = - \dfrac{\bd}{\rho} \cdot \dfrac{r_6 \left( r_1 \e^{r_1 x} - r_2 \e^{r_2 x} \right)}{(r_1 - r_6) \e^{r_1 x} - (r_2 - r_6) \e^{r_2 x}}+ \gam\left( \dfrac{\bd}{\rho}  \right)^2 \dfrac{r^2_6 (\e^{r_1 x} - \e^{r_2 x})( r_1 \e^{r_1 x} - r_2 \e^{r_2 x})}{\big((r_1 - r_6) \e^{r_1 x} - (r_2 - r_6) \e^{r_2 x}\big)^2} \\
			- \dfrac{\gam}{2} \left( \dfrac{\bd}{\rho} \right)^2 \dfrac{\big( r_8(r_1 + r_6) - 2r_1 r_6 \big) \e^{r_1 x} - \big( r_8(r_2 + r_6) - 2r_2 r_6 \big) \e^{r_2 x}}{(r_1 - r_6) \e^{r_1 x} - (r_2 - r_6) \e^{r_2 x}} \dfrac{r_3 \e^{r_3 x} - r_4 \e^{r_4 x}}{(r_3 - r_8) \e^{r_3 x} - (r_4 - r_8) \e^{r_4 x}}. 
		\end{align*}
		Assume the inequality in \eqref{eq:xt_pos} holds. Then, there exists a positive number $\eps$ such that for all $\gam \in (0, \eps)$, the equation $f(x, \gam) - 1 = 0$ admits a unique positive solution $\xt_\gam $ (that is, $f(\xt_\gam, \gam) - 1=0$ for some $\xt_\gam > 0$). 
	\end{lemma}
	
	\begin{proof} When $\gam = 0$, we easily see that $f(x, 0) - 1 = 0$ has a unique positive solution $\xt_0>0$ if (and only if)  \eqref{eq:xt_pos} holds.
		Next, by a tedious calculus, we verify that $\frac{\partial f}{\partial x}|_{(x, \gam) = (\xt_0, 0)} \neq 0$. Then, by the implicit function theorem, there exists a small positive $\eps > 0$ such that $f(x, \gam) - 1 =0$ has a unique solution $\xt_\gam>0$ for all $\gam \in (0, \eps)$.
	\end{proof}
	
\begin{figure}[htb]
		\centering
		\includegraphics[width = 0.6\textwidth, trim = 2cm 0cm 2cm 0cm, clip=true]{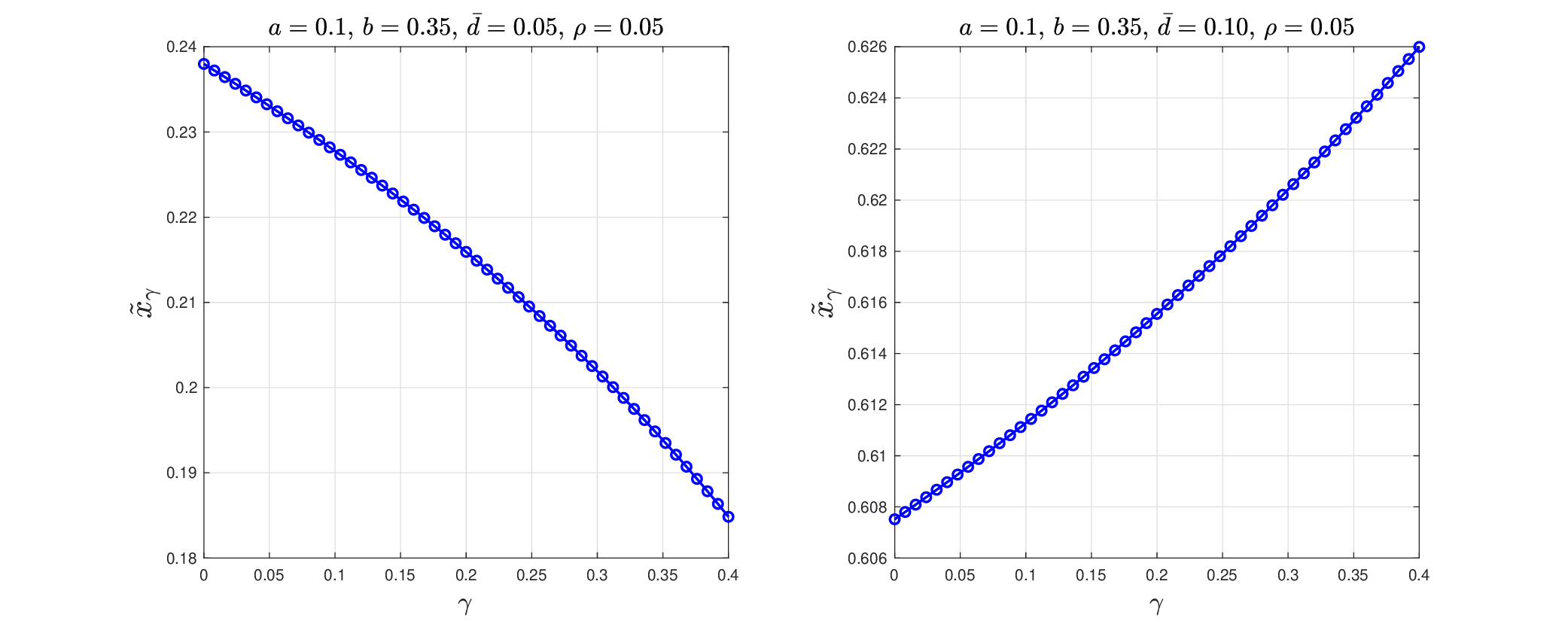}
		\includegraphics[width = 0.3\textwidth, trim = 2cm 0cm 1cm 0cm, clip=true]{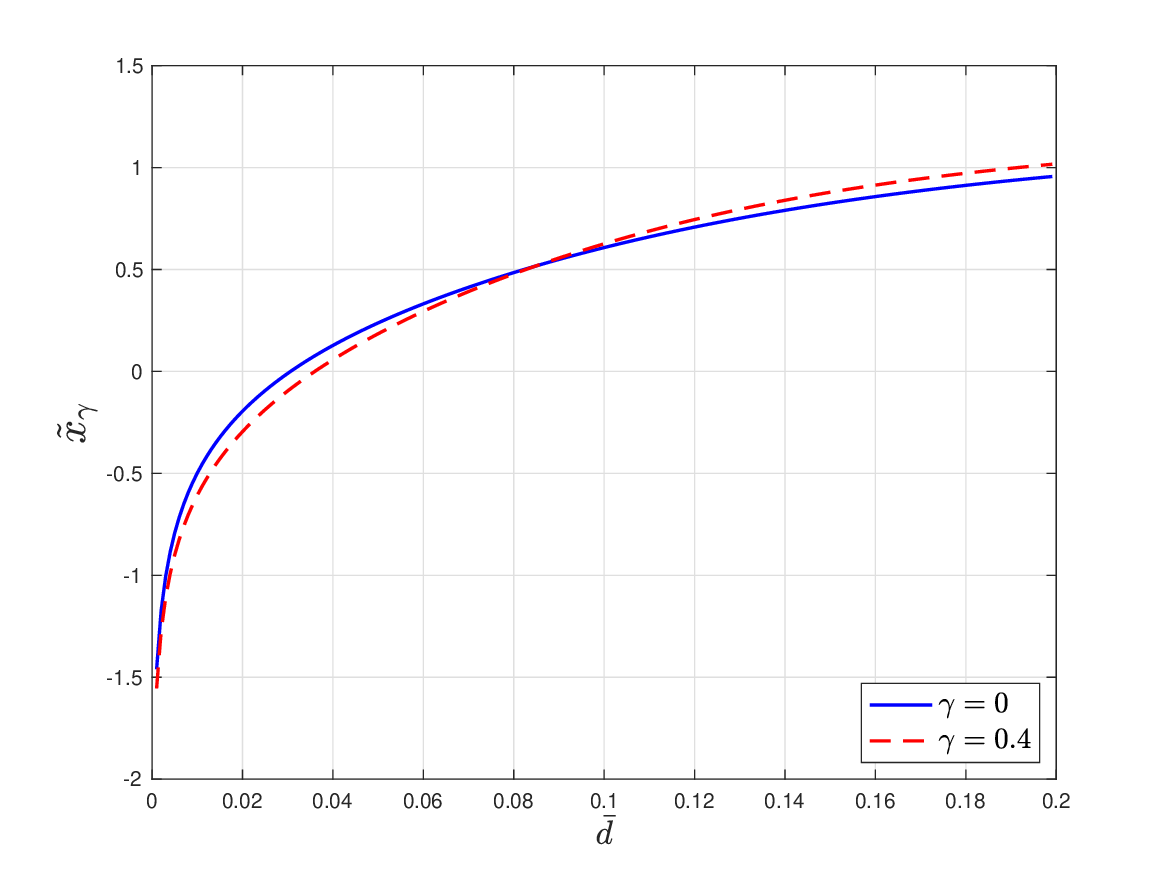}
		\\[-2ex]
		\caption{The unique solution $\xt_\gam$ to $f(x, \gam) - 1 = 0$ as a function of risk aversion $\gam$ or maximum rate $\bd$}
		\label{fig:gamma}
	\end{figure}

	Although \cref{lem:small} does not provide a precise bound on $\eps$, a standard nonlinear solver can easily determine whether $f(x, \gam) - 1 = 0$ has a unique (positive) solution once the model parameters are given. As an example, we set $a=0.1$, $b=0.35$, $\rho = 0.05$, and $\bar{d} = 0.05$ or $0.1$ and find that a unique solution $\xt_\gam$ exists for all $\gam \in [0, 0.4]$. We plot $\xt_\gam$ as a function of $\gam$ in the first two panels of \cref{fig:gamma} ($\bar{d} = 0.05$ in the left panel and $\bar{d} = 0.1$ in the middle panel). Notably, $\xt_\gam$ 
	exhibits monotonicity with respect to $\gam$, but its direction, whether increasing or decreasing, depends on the parameters.  
	In addition, we plot $\tilde{x}_\gam$ as a function of the maximum rate $\bd$ for two cases of $\gam$, $\gamma=0$ and $\gamma = 0.4$, in the right panel of \cref{fig:gamma}. We observe that $\tilde{x}_\gam$ increases with respect to $\bd$, which is consistent with our intuition.

	Motivated by the explicit result when $\gam = 0$, we construct a barrier strategy for small positive $\gam$ using the solution $\xt_\gam$ from \cref{lem:small} and show in the next theorem that it is the equilibrium strategy in the sense of \cref{def:bq_timeconsist}. 
	
	\begin{theorem}\label{thm:small_gam_bd}
		Assume that the equation $f(x, \gam) - 1 = 0$ has a unique positive solution $\xt := \xt_\gam > 0$ for the company's risk aversion $\gam$, as established in \cref{lem:small}. 
		Then, a time-independent dividend strategy induced by 
		\begin{equation}\label{eq:pdt}
			\pd^*(x) =
			\begin{cases}
				0, &\quad  x \le \xt, \\
				\bd, &\quad  x > \xt,
			\end{cases}
		\end{equation}
		is a time-consistent equilibrium dividend strategy in \cref{def:bq_timeconsist}. 
		
		Moreover, $G(x) = \E_{x}(Y_0^*)$ and $H(x) = \E_{x} ((Y_0^*)^2)$  are given by 
		\begin{align}
			\label{eq:G2}
			G(x) = C_1 \big(\e^{r_1 x} - \e^{r_2 x} \big) \cdot  \id_{ \{ x \le \xt \} } + \big(C_6 \e^{r_6 x} + \bd/\rho \big) \cdot  \id_{ \{ x > \xt \} },
		\end{align}
		in which $C_1$ and $C_6$ are given by \eqref{eq:C1}, 
		and 
		\begin{align}\label{eq:H2}
			H(x) = C_3 \big(\e^{r_3 x} - \e^{r_4 x} \big) \cdot  \id_{ \{ x \le \xt \} } + \left(	C_8 \e^{r_8 x} + 2 \, (\bd/\rho) \, C_6 \e^{r_6 x} + (\bd/\rho)^2  \right) \cdot  \id_{ \{ x > \xt \} },
		\end{align}
		in which $C_3$ and $C_8$ are given by \eqref{eq:C3} and \eqref{eq:C8}, respectively. 
		The corresponding value function equals $V(x, t) \equiv V(x, 0) = G(x) - \frac{\gam}{2} \left( H(x) - G^2(x) \right)$.
	\end{theorem}
	
	\begin{proof}
		Because the MV dividend problem in \cref{def:bq_timeconsist} is time-homogeneous, the value function $V(x, t)$ is time-independent; with slight abuse of notation, we write $V(x)$ for $V(x, t)$ and use $V'$ and $V''$ to denote its first and second derivatives. (The same rule applies to $G(x)$ and $H(x)$.) By examining the HJB equation in \eqref{eq:HJB}, we find that the sup problem is \emph{linear} with respect to $d$. Therefore, we hypothesize that a barrier strategy in the form of \eqref{eq:pdt} is an equilibrium dividend strategy, 
		with the barrier $\xt>0$ yet to be identified.
		
		With the candidate strategy in \eqref{eq:pdt}, we proceed to derive $G(x)$ and $H(x)$ from \eqref{eq:G} and \eqref{eq:H} in the extended HJB system. First, consider $x \le \xt$ and note from \eqref{eq:pdt} that $\pd^*(x) = 0$. In this case, $G$ satisfies the following boundary-value ODE: $- \rho G(x) + a G'(x) + \frac{1}{2} \, b^2 G''(x) = 0$, with $G(0) = 0$. The solution is given by the first part in \eqref{eq:G2}, with $C_1>0$ yet to be determined. By a similar argument, we obtain the solution of $H$ for $x \le \xt$ in \eqref{eq:H2}. Next, we study the case of $x > \xt$ and note that now $\pd^*(x) = \bd$. The ODE satisfied by $G$ in this case becomes $- \rho G(x) + (a - \bd) G'(x) + \frac{1}{2} \, b^2 G''(x) + \bd= 0$, and we obtain the general solution $G(x) = C_5 \e^{r_5 x} + C_6 \e^{r_6 x} + \frac{\bd}{\rho}$, in which $r_5$ and $r_6$ are defined by \eqref{eq:r_one}, and $C_5$ and $C_6$ are constants yet to be determined. Because $r_5 > 0$, the growth condition on $G$ in \cref{thm:verif_bd} requires that $C_5 = 0$, which implies the second part in \eqref{eq:G2}. By substituting $\pd^* = \bd$ and the solution of $G$ into \eqref{eq:H}, we solve the ODE of $H$ and verify the second part in \eqref{eq:H2}; note that the constant $C_6$ in $H$ is the same from $G$. 
		
		With $G$ and $H$ identified in \eqref{eq:G2} and \eqref{eq:H2}, respectively, we still need to deduce the constants $C_1$, $C_3$, $C_6$, and $C_8$, and we achieve that by ``smooth pasting.'' Recall that $G$ and $H$ are $\mC^2$, except they may only be $\mC^1$ at $\xt$. As such, we establish the smooth pasting conditions of $G$ and $H$ by (1)  $G(\xt^-) = G(\xt^+)$ and  $G'(\xt^-) = G'(\xt^+)$, 
		from which we obtain 
		{\small 
			\begin{align}
				\label{eq:C1}
				C_1 = - \, \dfrac{\bd}{\rho}   \dfrac{r_6}{(r_1 - r_6) \e^{r_1 \xt} - (r_2 - r_6) \e^{r_2 \xt}} > 0 \text{ and } C_6 = - \, \dfrac{\bd}{\rho} \dfrac{r_1 \e^{(r_1 - r_6)\xt} - r_2 \e^{(r_2 - r_6)\xt}}{(r_1 - r_6) \e^{r_1 \xt} - (r_2 - r_6) \e^{r_2 \xt}} < 0,	
			\end{align}
		}
		and (2) $H(\xt^-) = H(\xt^+)$ and $H'(\xt^-) = H'(\xt^+)$, from which we obtain 
		{\small 
			\begin{align}
				C_3 = \left( \dfrac{\bd}{\rho} \right)^2 \dfrac{\big( r_8(r_1 + r_6) - 2r_6 r_1 \big) \e^{r_1 \xt} - \big( r_8(r_2 + r_6) - 2r_6 r_2 \big) \e^{r_2 \xt}}{\big( (r_1 - r_6) \e^{r_1 \xt} - (r_2 - r_6) \e^{r_2 \xt} \big) \big( (r_3 - r_8) \e^{r_3 \xt} - (r_4 - r_8) \e^{r_4 \xt} \big)}, \label{eq:C3}
			\end{align}
		}
		and 	
		{\small 
			\begin{align}
				C_8 &= \left( \dfrac{\bd}{\rho} \right)^2 \dfrac{ \e^{(r_3 - r_8) \xt} - \e^{(r_4 - r_8) \xt}}{(r_3 - r_8) \e^{r_3 \xt} - (r_4 - r_8) \e^{r_4 \xt}} \cdot \dfrac{\big( r_8(r_1 + r_6) - 2r_6 r_1 \big) \e^{r_1 \xt} - \big( r_8(r_2 + r_6) - 2r_6 r_2 \big) \e^{r_2 \xt}}{(r_1 - r_6) \e^{r_1 \xt} - (r_2 - r_6) \e^{r_2 \xt}} \notag \\
				&\qquad + \left( \dfrac{\bd}{\rho} \right)^2 \dfrac{(r_1 + r_6) \e^{(r_1 - r_8) \xt} - (r_2 + r_6) \e^{(r_2 - r_8) \xt}}{(r_1 - r_6) \e^{r_1 \xt} - (r_2 - r_6) \e^{r_2 \xt}}. \label{eq:C8}
			\end{align}
		}

		Next, to characterize the barrier $\xt$, we use $V'(\xt) = 1$, which, after some tedious calculus, is equivalent to $f(\xt, \gam) - 1 = 0$, whose solution is established by \cref{lem:small}. Note that the condition $V'(\xt) = 1$ naturally implies $V \in \mC^2$ by \eqref{eq:HJB}.
		
		By construction, $V$, $G$, and $H$ satisfy all the conditions of  \cref{thm:verif_bd}, and the strategy induced by $\pd^*$ in \eqref{eq:pdt} is admissible. Therefore, by the verification theorem in \cref{thm:verif_bd}, all assertions of \cref{thm:small_gam_bd} follow.
	\end{proof}

	Recall from Taksar \cite{T2000} that when $\gam = 0$, if $\frac{\bd}{\rho} + \frac{1}{r_6} < 0$, it is always optimal to pay at the maximum rate $\bd$. Below, we show that this result extends to small positive $\gam$. Recall that all the $r_i$s are defined in \eqref{eq:r_one} and \eqref{eq:r_two}.

	\begin{theorem}
		\label{thm:max_d}
		Assume $\frac{\bd}{\rho} + \frac{1}{r_6} < 0$. 
		There exists $\tilde{\eps} > 0$ such that, if $\gam < \tilde{\eps}$,
		then paying dividends at the maximum rate $\bd$ is a time-consistent equilibrium strategy (that is, $\pd^*(x) = \bd$ for all $x > 0$).   
		Moreover, $G(x) = \E_{x}(Y_0^*)$ and $H(x) = \E_{x} ((Y_0^*)^2)$  are given by
		\begin{align}
			\label{eq:GH_max}
			G(x) =  \left( 1 - \e^{r_6 x} \right) \, \bd/\rho \quad \text{ and } \quad 
			H(x)  =  \left(1 - 2 \e^{r_6 x} + \e^{r_8x} \right) \, (\bd/\rho)^2,
		\end{align}
		and the equilibrium value function equals $	V(x) 
		= G(x) - \frac{\gam}{2} (H(x) - G^2(x)).$
	\end{theorem}
	
	\begin{proof}
		Assume for a moment that $\frac{\bd}{\rho} + \frac{1}{r_6} < 0$ is a sufficient condition for $\pd^*(x) = \bd$ for all $x > 0$. Under this hypothesis, $G$ solves the ODE of $- \rho G(x) + (a - \bd) G'(x) + \frac{1}{2} \, b^2 G''(x) + \bd = 0$ over $x \in (0, \infty)$, with the boundary condition $G(0) = 0$; solving this ODE, along with the growth condition in \cref{thm:verif_bd}, leads to the unique solution of $G$ in \eqref{eq:GH_max}. Then, $H$ solves $- 2 \rho H(x) + (a - \bd) H'(x) + \frac{1}{2} \, b^2 H''(x) = - \, \frac{2\bd^2}{\rho} \, \big( 1 - \e^{r_6 x} \big)$, with $H(0) = 0$; by a similar argument, we verify that the solution is uniquely given by the expression of $H$ in \eqref{eq:GH_max}. 
		
		Using $G$ and $H$ in \eqref{eq:GH_max}, we immediately get the candidate value function $V$, which satisfies the HJB equation in  \eqref{eq:HJB} if the optimizer is $\bd$. As such, the remaining task is to show that the supremum in \eqref{eq:HJB} is indeed achieved at $\bd$, which is equivalent to $V'(x) \le 1$ for all $x > 0$.  For that purpose, note that the variance term, $\V_x(Y_0^*) = \left(\bd/\rho\right)^2 \left(\e^{r_8x} - \e^{2 r_6 x}\right)$, first increases to its maximum at $x_m = \frac{1}{r_8-2r_6} \ln \frac{2r_6}{r_8} > 0$ and eventually decreases to 0. In addition, this term is concave for $x < 2 x_m$ and convex for $x > 2 x_m$; thus, $V$ is concave ($V'' (x) \le 0$) for $x \ge 2 x_m$, which means we only need to show that $V'(x) \le 1$ for all $x <2x_m$. We compute $V'(x) = - \, \frac{\bd}{\rho} \, r_6 \e^{r_6 x} - \frac{\gam}{2} \left( \frac{\bd}{\rho} \right)^2 \big(r_8 \e^{r_8x} - 2r_6 \e^{2r_6x} \big)$.
		It is easy to see that the second term of $V'(x)$, excluding the factor $- \frac{\gam}{2}$, is bounded on $(0, 2 x_m)$.
		If the assumption $\frac{\bd}{\rho} + \frac{1}{r_6} < 0$ holds, the first term of $V'(x)$ is strictly less than 1, by recalling $r_6 < 0$ from \eqref{eq:r_one}. Therefore, there exists $\tilde{\eps} > 0$ such that $V'(x) \le 1$  over $(0, 2 x_m)$ for all $\gam < \tilde{\eps}$, which, along with $V''(x) \le 0$ for $x \ge 2 x_m$, proves $\pd^*(x) = \pd$ for all $x > 0$.  The proof is now completed.
	\end{proof}
	
	For the same model parameters in Figure \ref{fig:gamma}, $\frac{\bd}{\rho} + \frac{1}{r_6} < 0$ is equivalent to $\bar{d}<0.0306$. By setting $\bar{d} = 0.03$, we numerically find that $\tilde{\eps} \approx 0.65$, and for all $\gamma \in [0, \tilde{\eps})$, the solution to $f(x, \gam) - 1 = 0$, $\tilde{x}_\gam$, is strictly negative, and it decreases with respect to $\gam$ (computation results are available upon request). However, the results in Theorem \ref{thm:max_d} might not hold for sufficiently large $\gam > 0$, differing from that when $\gam = 0$ in Taksar \cite{T2000}.

	\section{Conclusion}
	\label{sec:con}
	
	This paper formulates an optimal dividend rate problem under the MV criterion for a company whose (uncontrolled) surplus is modeled by a Brownian motion with drift. Due to the time-inconsistency arising from the MV objective, we seek a time-consistent equilibrium dividend strategy. We prove a new verification theorem that is tailored to the MV dividend problem and apply it to show that for small $\gam$ (risk aversion toward variance), a barrier dividend strategy is the equilibrium strategy. 
	
	Two successful applications of the verification lemma, both yielding close-form solutions, require a small risk aversion $\gam$ (see Theorems \ref{thm:small_gam_bd} and \ref{thm:max_d}).  
	It remains an open question to derive the equilibrium strategy for all $\gam > 0$. Here, the ``penalty'' on the dividend variability is via the variance term in the optimization objective; it could be interesting to explore the possibility of directly penalizing the company's surplus when it follows a volatile dividend strategy.
	Furthermore, we adopt the classical control framework in this paper, and the company's dividend rates are bounded. In future research, we will explore the singular control setup and allow unbounded dividend rates (lump-sum payments).

	\section*{Acknowledgments}
	We thank two anonymous reviewers for their comments.



\begin{thebibliography}{50}
		
		\bibitem{AT2009} Albrecher, Hansj\"org and Stefan Thonhauser (2009). Optimality results for dividend problems in insurance. {\it Revista de la Real Academia de Ciencias Exactas, Fisicas y Naturales}, 103(2), 295-320.
		
		\bibitem{AAM2022} Albrecher, Hansj\"org, Pablo Azcue, and Nora Muler (2022). Optimal ratcheting of dividends in a Brownian risk model. {\it SIAM Journal on Financial Mathematics}, 13(3), 657-701.
		
		\bibitem{ABY2019} Angoshtari, Bahman, Erhan Bayraktar, and Virginia R. Young (2019). Optimal dividend distribution under drawdown and ratcheting constraints on dividend rates. {\it SIAM Journal on Financial Mathematics}, 10(2), 547-577.
		
		\bibitem{AT1997} Asmussen, Soren and Michael Taksar (1997). Controlled diffusion models for optimal dividend pay-out. {\it Insurance: Mathematics and Economics}, 20(1), 1-15.
		
		\bibitem{A2009} Avanzi, Benjamin (2009). Strategies for dividend distribution: A review. {\it North American Actuarial Journal}, 13(2), 217-251.
		
		\bibitem{AFS2023} Avanzi, Benjamin, Debbie Kusch Falden, and Mogens Steffensen (2023). Stable dividends under linear-quadratic optimisation. {\it Quantitative Finance}, 23(9), 1199-1215.
		
		\bibitem{B2009} Baker, H. Kent (2009). {\it Dividends and Dividend Policy}, Volume 1. John Wiley \& Sons.
		
		\bibitem{BC2010} Basak, Suleyman and Georgy Chabakauri (2010). Dynamic mean-variance asset allocation. {\it Review of Financial Studies}, 23(8), 2970-3016.
		
		\bibitem{BM2010} Bj\"ork, Tomas and Agatha Murgoci (2010). A general theory of Markovian time inconsistent stochastic control problems. Working paper, available at SSRN 1694759.
		
		\bibitem{BMZ2014} Bj\"ork, Tomas, Agatha Murgoci, and Xun Yu Zhou (2014). Mean-variance portfolio optimization with state-dependent risk aversion. {\it Mathematical Finance}, 24(1), 1-24.
		
		\bibitem{CLZ2014} Chen, Shumin, Zhongfei Li, and Yan Zeng (2014). Optimal dividend strategies with time-inconsistent preferences. {\it Journal of Economic Dynamics and Control}, 46, 150-172.
		
		\bibitem{CLZ2018} Chen, Shumin, Zhongfei Li, and Yan Zeng  (2018). Optimal dividend strategy for a general diffusion process with time-inconsistent preferences and ruin penalty. {\it SIAM Journal on Financial Mathematics}, 9(1), 274-314.
		
		\bibitem{D1957} De Finetti, Bruno. Su un’impostazione alternativa della teoria collettiva del rischio. {\it Transactions of the XVth International Congress of Actuaries}, 2(1), 433-443.
		
		\bibitem{FS2006} Fleming, Wendell H. and H. Mete Soner (2006). {\it Controlled Markov Processes and Viscosity Solutions}, second edition, Springer.
		
		
		\bibitem{G1991} Grandell, Jan (1991).  {\it Aspects of Risk Theory}. Springer-Verlag, New York.
		
		\bibitem{GX2024} Guan, Chonghu and Zuo Quan Xu (2024). Optimal ratcheting of dividend payout under Brownian motion surplus. {\it SIAM Journal on Control and Optimization}, 62(5), 2590-2620.
		
		\bibitem{JS1995} Jeanblanc-Picqu\'e, Monique and Albert N. Shiryaev (1995). Optimization of the flow of dividends. {\it Russian Mathematical Surveys}, 50(2), 257-278.
		
		\bibitem{KS2015} Kronborg, Morten Tolver, and Mogens Steffensen (2015). Inconsistent investment and consumption problems. {\it Applied Mathematics and Optimization}, 71, 473-515.
		
		\bibitem{LLLY2018} Landriault, David, Bin Li, Danping Li, and Virginia R. Young (2018). Equilibrium strategies for the mean-variance investment problem over a random horizon. {\it SIAM Journal on Financial Mathematics}, 9(3), 1046-1073.
		
		\bibitem{T2000} Taksar, Michael I. (2000). Optimal risk and dividend distribution control models for an insurance company. {\it Mathematical Methods of Operations Research}, 51, 1-42.
		
		\bibitem{ZSY2020} Zhu, Jinxia, Tak Kuen Siu, and Hailiang Yang (2020). Singular dividend optimization for a linear diffusion model with time-inconsistent preferences. {\it European Journal of Operational Research}, 285(1), 66-80.
		
		\bibitem{ZJ2020} Zhou, Zhou, and Zhuo Jin (2020). Optimal equilibrium barrier strategies for time-inconsistent dividend problems in discrete time. {\it Insurance: Mathematics and Economics}, 94, 100-108.
		
	\end{thebibliography}
	
\end{document}